\documentclass[a4paper,11pt]{article}

\usepackage[left=1.5cm,right=1.5cm,top=1.5cm]{geometry}

\usepackage{graphicx}
\usepackage[colorlinks=true,citecolor=blue]{hyperref}

\usepackage{amsmath}
\usepackage{amsfonts}
\usepackage{amssymb}
\usepackage{color}
\usepackage{accents}

\usepackage{ssgc}

\usepackage{natbib}

\allowdisplaybreaks[2]

\begin{document}

\title{Granger causality for state space models}

\author{Lionel Barnett and Anil K. Seth\\Sackler Centre for Consciousness Science\\School of Engineering and Informatics\\University of Sussex, BN1 9QJ, UK}

\date{\today}

\maketitle

\begin{abstract}
Granger causality, a popular method for determining causal influence between stochastic processes, is most commonly estimated via linear autoregressive modeling. However, this approach has a serious drawback: if the process being modeled has a moving average component, then the autoregressive model order is theoretically infinite, and in finite sample large empirical model orders may be necessary, resulting in weak Granger-causal inference. This is particularly relevant when the process has been filtered, downsampled, or observed with (additive) noise - all of which induce a moving average component and are commonplace in application domains as diverse as econometrics and the neurosciences. By contrast, the class of autoregressive moving average models---or, equivalently, linear state space models---is closed under digital filtering, downsampling (and other forms of aggregation) as well as additive observational noise. Here, we show how Granger causality, conditional and unconditional, in both time and frequency domains, may be calculated simply and directly from state space model parameters, via solution of a discrete algebraic Riccati equation. Numerical simulations  demonstrate that Granger causality estimators thus derived have greater statistical power and smaller bias than pure autoregressive estimators. We conclude that the state space approach should be the default for (linear) Granger causality estimation.
\end{abstract}

\section{Introduction} \label{sec:intro}

Wiener-Granger causality (GC), a powerful method for determining information flow or directed functional connectivity between stochastic variables, is based on the premise that cause (a) precedes effect, and (b) contains unique information about effect \citep{Wiener:1956,Granger:1963,Granger:1969}. Since its inception, it has steadily gained popularity in a broadening range of fields (see \eg\ \citealt{SethJNEURO:2015}), due to its data-driven nature (few structural assumptions need be made about the data generation process), conceptual simplicity, spectral decomposition property and ease of implementation. It is most commonly operationalized in a multivariate linear autoregressive (AR) modeling context, to the extent that it is frequently viewed as an essentially linear autoregressive method. This is unfortunate, not just because it obscures the nonparametric essence of the original idea (indeed GC may be formalized in purely information-theoretic terms \citep{HlavackovaEtal:2007,AmblardMichel:2011,Barnett:teml:2012}), but because the pure AR modeling approach frequently sits uncomfortably with the structure of the time series data to which it is applied.

Time series data from diverse application domains, in particular econometrics and the neurosciences, often contain a strong moving average (MA) component, which may not be represented parsimoniously by a finite order AR model. Apart from any MA component intrinsic to the underlying signal, an MA component may be induced in an observed process by common data acquisition, sampling and pre-processing procedures. Thus, for example, if a (finite order) AR process is subsampled or aggregated, or observed with additive measurement noise, the resultant process will be an autoregressive moving average (ARMA) process \citep{Solo:2007}. A subprocess of an AR process will also generally be ARMA \citep{NsiriRoy:1993}, as will a filtered AR process \citep{Barnett:gcfilt:2011}. In general an ARMA process will have \emph{infinite} AR model order - but of course in finite sample a finite, possibly large, model order must be selected, which will be reflected in GC statistics of reduced statistical power and increased bias.

In contrast to linear AR models, the class of multivariate ARMA models---or, equivalently \citep{HandD:2012}, finite order linear \emph{state space} (SS) models---is closed under all of the above-mentioned operations. While the potential of SS modeling for GC inference has been remarked on \citep{Solo:2007,ValdesSosaEtal:2011,Seth:gcfmri:2013}, a rigorous derivation and demonstration of this potential has so far been lacking. Here we show how GC, conditional and unconditional, in both time and frequency domains, may be easily derived from SS parameters via solution of a single discrete algebraic Riccati equation (DARE). We verify in simulation, for a minimal ARMA process, the increase in statistical power and reduction in bias of GC estimated via SS, as compared to AR, modeling. We also discuss potential extensions of the SS approach beyond the usual stationary linear scenario, to encompass cointegrated processes, models with time-varying parameters, and nonlinear models. State space Granger causal inference stands to significantly enhance our ability to identify and understand information flow and causal interactions in a wide range of complex dynamical systems. Given the ready availability of efficient and effective state space system identification procedures, state space modeling should become the default approach to Granger causal analysis.

\section{State space models}

Our starting point\footnote{Notational conventions: vector quantities are written in lower-case bold, matrices in upper case. Superscript `$\trop$' denotes the transpose, `$*$' conjugate transpose and $|\cdot|$ the determinant of a (complex) matrix; superscript `$\rrop$' refers to a ``reduced'' model. $\expect{\,\cdot}$ denotes expectation and $\cexpect{\,\cdot}{\cdot}$ conditional expectation.} is a a discrete-time, real-valued  vector stochastic process $\by_t = \bracs{y_{1t} \ y_{2t} \ \ldots \ y_{nt}}^\trop$, $-\infty < t < \infty$, of observations. The general time-invariant linear SS model without input for the observation process $\by_t$ is
\begin{subequations}
\begin{align}
	\bx_{t+1} &= A\bx_t + \bu_t && \text{state transition equation} \label{eq:ssps} \\
	\by_t     &= C\bx_t + \bv_t && \text{observation equation} \label{eq:sspo}
\end{align} \label{eq:ssp}%
\end{subequations}
where $\bx_t$ is an (unobserved) $m$-dimensional state variable, $\bu_t,\bv_t$ are zero-mean white noise processes, $C$ is the observation matrix and $A$ the state transition matrix. The parameters of the model \eqref{eq:ssp} are $(A,C,Q,R,S)$, where
\begin{equation}
	\begin{bmatrix} Q & S \\ \tS & R  \end{bmatrix} \equiv
	\expect{ \hspace{-1pt}
	\begin{bmatrix} \bu_t \\ \bv_t \end{bmatrix} \hspace{-6pt}
 	\begin{array}{c}
		\begin{bmatrix} \bu_t^\trop & \!\! \bv_t^\trop \end{bmatrix} \\ \phantom.
	\end{array} \hspace{-4pt}
	}
\end{equation}
is the noise covariance matrix. We assume that $\bx_t,\by_t$ are weakly stationary, which requires that the transition equation \eqref{eq:ssps} be \emph{stable}: the stability condition is $\lambda_{\text{max}}(A) < 1$, where $\lambda_{\text{max}}(A)$ denotes the maximum of the absolute values of the eigenvalues of $A$. We also assume that $R$ is positive-definite. A process $\by_t$ satisfying the stable SS model \eqref{eq:ssp} also satisfies a stable ARMA model; conversely, any stable ARMA process may be shown to satisfy a stable SS model of the form \eqref{eq:ssp} \citep{HandD:2012}.

Given an SS model in general form \eqref{eq:ssp}, we define $\bz_t \equiv \cexpect{\bx_t}{\by^-_{t-1}}$, the projection of the state variable $\bx_t$ on the space spanned by the infinite past $\by^-_{t-1} \equiv \bracs{\by^\trop_{t-1} \ \by^\trop_{t-2} \ \ldots}^\trop$ of the observation variable \citep{HandD:2012}. It is then not difficult to show \citep{Aoki:1994} that the \emph{innovations} $\beps_t \equiv \by_t - \cexpect{\by_t}{\by^-_{t-1}}$ constitute a white noise process with positive-definite covariance matrix $\Sigma \equiv \expect{\beps^{\phantom\trop}_t \beps^\trop_t}$, and that in terms of the new state variable $\bz_t$ we have an SS model
\begin{subequations}
\begin{align}
	\bz_{t+1} &= A\bz_t + K \beps_t && \text{state transition equation} \label{eq:ssis} \\
	\by_t     &= C\bz_t + \beps_t && \text{observation equation} \label{eq:ssio}
\end{align} \label{eq:ssi}%
\end{subequations}
for $\by_t$, where $K$ is the Kalman gain matrix. The SS model \eqref{eq:ssi} is said to be in \emph{innovations form}, with parameters $(A,C,K,\Sigma)$. Note that innovations form constitutes a special case of \eqref{eq:ssp} with $\bu_t = K\beps_t, \bv_t = \beps_t$, so that $Q = K \Sigma \tK$, $ S = K \Sigma$ and $R = \Sigma$.

We can write the state equation \eqref{eq:ssis} as $\bz_t = (I - A \lag)^{-1} K \lag \cdot \beps_t$ where here $\lag$ represents the back-shift operator\footnote{In the frequency domain $z = e^{-i\omega}$ where $-\pi < \omega \le \pi$ is the phase angle (we note that the inverse $\lag^{-1}$ is sometimes used for the back-shift operator, particularly in the signal processing literature).}, which yields the MA representation for the observation process
\begin{equation}
	\by_t = H(\lag) \cdot \beps_t \,, \quad H(z) \equiv I + C (I - A z)^{-1} K z \label{eq:ssioma}
\end{equation}
with transfer function $H(z)$. From \eqref{eq:ssi} we have $\bz_{t+1} = \AKC \bz_t + K \by_t$ with $\AKC \equiv A-KC$, from which we may derive the AR representation
\begin{equation}
	\AKC(\lag) \cdot \by_t = \beps_t \,, \quad \AKC(z) \equiv I - C (I - \AKC z)^{-1} K z \label{eq:ssioar}
\end{equation}
where $\AKC(z) = H(z)^{-1}$ is the inverse transfer function. The model is \emph{minimum phase} if $\lambda_{\text{max}}(\AKC) < 1$ (a stable inverse ARMA model for the process $\by_t$ then exists); we assume minimum phase from now on. From \eqref{eq:ssioma} the cross-power spectral density (CPSD) of the observation variable has the factorization \citep{Wilson:1972}
\begin{equation}
	S(z) = H(z) \Sigma H^*(z) \label{eq:cpsdi}
\end{equation}
on $|z| = 1$ in the complex plane.

Assuming stability, minimum phase and positive-definite observation noise covariance\footnote{In fact these conditions may be relaxed somewhat - see \eg\ \citep{ChanEtal:1984,Solo:2015}.}, the DARE \citep{LancasterRodman:1995}
\begin{equation}
	P - A P \tA = Q - (A P \tC + S) (C P \tC + R)^{-1} (C P \tA + \tS) \label{eq:Pdare}
\end{equation}
has a  unique stabilizing solution for $P$, and we have \citep{Kailath:1980}
\begin{subequations}
\begin{align}
	\Sigma &= C P \tC + R \label{eq:p2iSig} \\
	K &= (A P \tC +S) \Sigma^{-1} \label{eq:p2iK}
\end{align} \label{eq:p2i}%
\end{subequations}
Given general SS parameters $(A,C,Q,R,S)$ then, corresponding innovations form parameters $K,\Sigma$ may be obtained through \eqref{eq:p2i} via solution of \eqref{eq:Pdare}.

\section{Granger causality} \label{sec:gc}

Granger causality is commonly expressed in terms of \emph{prediction error}. The full multivariate, conditional and spectral theory of GC in an AR framework was consolidated by \citet{Geweke:1982,Geweke:1984}, whose formulation we follow here.

Suppose that an observable process $\by_t$ is partitioned into subprocesses: $\by_t = \bracs{\by^\trop_{1t} \ \by^\trop_{2t} \ \by^\trop_{3t}}^\trop$. Within the ``observable universe of information'' represented by the process $\by_t$, GC from the $\by_{2t}$ to $\by_{1t}$ (conditional on $\by_{3t}$) quantifies the extent to which the past of $\by_{2t}$ improves prediction of the future of $\by_{1t}$ over and above the extent to which $\by_{1t}$ (along with $\by_{3t}$) already predicts its own future. Now the best prediction (in the least-squares sense) of $\by_{1t}$, given the entire universe of past information $\by^-_{t-1}$, is the projection $\cexpect{\by_{1t}}{\by^-_{t-1}}$. This prediction may be contrasted with the prediction $\cexpect{\by_{1t}}{\by^{\rrop-}_{t-1}}$ of $\by_{1t}$ based on the \emph{reduced} universe of past information $\by^{\rrop-}_{t-1}$, where $\by^\rrop_t \equiv \bracs{\by^\trop_{1t} \ \by^\trop_{3t}}^\trop$ omits $\by_{2t}$ from the observable information set. In Geweke's formulation, predictive power is quantified by the \emph{generalized variances} \citep{Wilks:1932,Barrett:2010} of the associated full and reduced residual errors (innovations), $\beps_t \equiv \by_t - \cexpect{\by_t}{\by^-_{t-1}}$ and $\beps^\rrop_t \equiv \by^\rrop_t - \cexpect{\by^\rrop_t}{\by^{\rrop-}_{t-1}}$ respectively: specifically, \citet{Geweke:1984} defines the time-domain GC from $\by_2$ to $\by_1$ conditional on $\by_3$ (for \emph{un}conditional GC we may take $\by_{3t}$ to be empty) as
\begin{equation}
	\cgc{\by_2}{\by_1}{\by_3} \equiv \ln\frac{|\Sigma^\rrop_{11}|}{|\Sigma_{11}|} \label{eq:cgctd}
\end{equation}
where $\Sigma = \expect{\beps^{\phantom\trop}_t {\beps^\trop_t}}$ and $\Sigma^\rrop = \expect{\beps^{\rrop\phantom\trop}_t\!\!\beps^{\rrop\trop}_t}$. In a maximum likelihood (ML) framework, the corresponding statistic is just the log-likelihood ratio \citep{NeymanPearson:1928,NeymanPearson:1933} for the nested AR models associated with the projections $\cexpect{\by_{1t}}{\by^-_{t-1}}$ and $\cexpect{\by_{1t}}{\by^{\rrop-}_{t-1}}$. It also has a natural interpretation as the rate of ``information transfer'' from process $\by_{2t}$ to process $\by_{1t}$ \citep{Palus:2001,Barnett:tegc:2009,Barnett:teml:2012}.

\citet{Geweke:1982,Geweke:1984} goes on to define GC $\csgc{\by_2}{\by_1}{\by_3}(z)$ in the spectral domain, and establishes the fundamental frequency decomposition
\begin{equation}
	 \cgc{\by_2}{\by_1}{\by_3} = \frac1{2\pi} \int_{-\pi}^\pi \csgc{\by_2}{\by_1}{\by_3}(e^{-i\omega}) \,d\omega \label{eq:gcspdecomp}
\end{equation}
In the unconditional case (\ie\ $\by_{3t}$ is empty) $\sgc{\by_2}{\by_1}(z)$ is given by \citep{Geweke:1982}
\begin{equation}
    \sgc{\by_2}{\by_1}(z) \equiv \ln\frac{|S_{11}(z)|}{|S_{11}(z) - H_{12}(z) \Sigma_{22|1} H^*_{12}(z)|} \label{eq:ugcfd}
\end{equation}
where $\Sigma_{ij|k} \equiv \Sigma_{ij} - \Sigma_{ik} \Sigma_{kk}^{-1} \Sigma_{kj}$ denotes a partial covariance matrix. The conditional case is somewhat more intricate: defining the process $\dby_t \equiv \bracs{\beps^{\rrop\trop}_{1t} \ \by^\trop_{2t} \ \beps^{\rrop\trop}_{3t}}^\trop$, the frequency-domain causality from $\by_2$ to $\by_1$ conditional on $\by_3$ is defined \citep{Geweke:1984} as the \emph{un}conditional causality
\begin{equation}
    \csgc{\by_2}{\by_1}{\by_3}(z) \equiv \sgc{\bracs{ \dby^\trop_2 \ \dby^\trop_3}^\trop}{\dby_1}(z) \label{eq:cgcfd}
\end{equation}
Let $\dS(z)$, $\dH(z)$ and $\dSigma$ be respectively the CPSD, transfer function and innovations covariance matrix of the process $\dby_t$. We note firstly that, since the innovations process $\beps^\rrop_t$ is white, $\dS_{11}(z)$ is just the flat spectrum $\Sigma^\rrop_{11}$. The AR representation for the reduced process $\by^\rrop_t$ is $B^\rrop(\lag) \cdot \by^\rrop_t = \beps^\rrop_t$ where $B^\rrop(z) \equiv {H^\rrop(z)}^{-1}$ is the inverse transfer function of the reduced model, so we have $\dby_t = \dB(\lag) \cdot \by_t$, where
\begin{equation}
    \dB(z) \equiv
	\begin{bmatrix}
		B^\rrop_{11}(z) & 0 & B^\rrop_{13}(z) \\
		0 & I & 0 \\
		B^\rrop_{31}(z) & 0 & B^\rrop_{33}(z) \\
	\end{bmatrix}
\end{equation}
But $\by_t = H(\lag) \cdot \beps_t$, so that $\dby_t = \dB(\lag) H(\lag) \cdot \beps_t$ and it follows that $\dH(z) = \dB(z) H(z)$ and $\dSigma = \Sigma$. From \eqreff{eq:ugcfd}{eq:cgcfd}, after some matrix algebra, we derive the following explicit expression for conditional GC in the frequency domain:
\begin{equation}
    \csgc{\by_2}{\by_1}{\by_3}(z) =
	\ln\frac{|\Sigma^\rrop_{11}|}{\left|
	\hspace{-8pt}
	\begin{array}{c}
	\vspace{-11pt} \\
	\begin{array}{c}
		\Sigma^\rrop_{11} -
		\begin{bmatrix} \dH_{12}(z) & \dH_{13}(z) \end{bmatrix} \\ \phantom.
	\end{array} \hspace{-6pt}
	\begin{bmatrix} \Sigma_{22|1} & \Sigma_{23|1} \\ \Sigma_{32|1} & \Sigma_{33|1} \end{bmatrix}
	\hspace{-5pt}
	\begin{bmatrix} \dH^*_{12}(z) \\ \dH^*_{13}(z) \end{bmatrix}
	\end{array}
	\hspace{-5pt}
	\right|} \label{eq:cgcfd1}
\end{equation}
where
\begin{equation}
 	\begin{array}{c}
		\begin{bmatrix} \dH_{12}(z) & \dH_{13}(z) \end{bmatrix} \equiv
		\begin{bmatrix} B^\rrop_{11}(z) & B^\rrop_{13}(z) \end{bmatrix} \\ \phantom.
	\end{array} \hspace{-8pt}
	\begin{bmatrix} H_{12}(z) & H_{13}(z) \\ H_{32}(z) & H_{33}(z) \end{bmatrix}
\end{equation}

Suppose now that we have an SS model in innovations form \eqref{eq:ssi} for the observation process $\by_t$, satisfying the assumptions of stability, minimum phase and positive-definite innovations covariance. The reduced process $\by^\rrop_t$ satisfies the observation equation:
\begin{equation}
	\by^\rrop_t = C^\rrop\bz_t + \bv^\rrop_t \label{eq:sspor}
\end{equation}
where $C^\rrop \equiv \bracs{C^\trop_1 \ C^\trop_3}^\trop$ and $\bv^\rrop_t \equiv \bracs{\beps^\trop_{1t} \ \beps^\trop_{3t}}^\trop$. Together with the original state transition equation \eqref{eq:ssis}, \eqref{eq:sspor} constitutes an SS model for $\by^\rrop_t$. Note that in general this ``reduced SS model'' will \emph{not} be in innovations form (it is for this reason that we write $\bv^\rrop_t$ rather than $\beps^\rrop_t$). It is, however, trivially stable (since $A$ is stable) and the noise covariance $R^\rrop \equiv \expect{\bv^{\rrop\phantom\trop}_t\!\!\bv^{\rrop\trop}_t}$ is positive-definite (since $\Sigma$ is positive-definite). By the minimum phase assumption, the full process $\by_t$ is invertible ARMA, so that the subprocess $\by^\rrop_t$ is also invertible ARMA \citep{NsiriRoy:1993} and the reduced SS model is thus minimum phase. Now the reduced Kalman gain $K^\rrop$ [which enters into the expression for $\AKC^\rrop(z)$], and the reduced innovations covariance $\Sigma^\rrop$, may be obtained by solving the DARE \eqref{eq:Pdare} for the general form SS \eqreff{eq:ssis}{eq:sspor}. Thus, given innovations form SS parameters $(A,C,K,\Sigma)$, GCs in both time and frequency domain, conditional and unconditional, are readily calculated.

$\cgc{\by_2}{\by_1}{\by_3}$ vanishes precisely when coefficients with block indices $12$ vanish at all lags in the AR representation of the process $\by_t$. For an SS model, setting $\AKC_{12}(\lag) \equiv 0$ in \eqref{eq:ssioar} yields the necessary and sufficient condition
\begin{equation}
	C_1 \AKC^k K_2 = 0 \,, \qquad k = 0,1,2,\ldots \label{eq:cgcnull}
\end{equation}
By the Cayley-Hamilton Theorem, \eqref{eq:cgcnull} need be satisfied just for $k = 0,1,\ldots,m-1$. In particular, $(B,C_1)$ \emph{observable} $\implies K_2 = 0$ while $(B,K_2)$ \emph{controllable} $\implies C_1 = 0$ (the latter case is trivial, since then $\by_{1t}$ is a pure white noise process uncorrelated with the pasts of $\by_{2t}$ and $\by_{3t}$). In terms of the MA representation \eqref{eq:ssioma}, performing block-inversion of $\AKC(\lag) = H(\lag)^{-1}$, the non-causality condition is found to be $H_{12}(\lag) -  H_{13}(\lag)  H_{33}(\lag)^{-1}  H_{32}(\lag) \equiv 0$ in the conditional case, or just $H_{12}(\lag) \equiv 0$ in the unconditional case (in which case \eqref{eq:cgcnull} holds with $\AKC$ replaced by $A$).

AR model parameters are generally estimated from data by standard regression techniques such as Ordinary Least Squares (OLS) or so-called LWR (Levinson-Wiggins-Robinson) algorithms \citep{Levinson:1947,Whittle:1963,WigginsRobinson:1965,MorfEtal:1978}, which yield pseudo-ML estimates \citep{Lutkepohl:2005}. In this study we are not primarily concerned with identification of SS models, for which a large literature exists \citep{Kailath:1980,Ljung:1999}. We note nonetheless that many identification procedures yield (asymptotically) pseudo-ML estimates for model parameters. \emph{State space-subspace} (SS-SS) algorithms, in particular, do not require iteration and are highly efficient \citep{VOandDM:1996}. See \apxref{sec:statinf} for further discussion on GC estimation and statistical inference for AR and SS models; the key point is that SS-SS GC estimation is in general more computationally efficient and numerically stable than AR GC estimation.

\section{Simulation experiment} \label{sec:simulation}

To examine the performance of state space Granger-causal inference, we generated and analyzed simulated time series data using a minimal AR process for which GC values can be computed analytically. The AR data was then filtered to induce a MA component. We compare the \emph{statistical power} and \emph{bias} of SS-derived and AR-derived Granger-causal estimators (see \apxref{sec:statpow}) for the resulting ARMA process.

To generate simulated time series, we used the bivariate AR($1$) process $\cby_t$ defined by \citep{Barnett:gcfilt:2011}
\begin{subequations} \label{eq:mvar}
\begin{align}
	\cy_{1t} &= a \cy_{1,t-1} + c \cy_{2,t-1} + \veps_{1t} \label{eq:minvar1} \\
	\cy_{2t} &= \hspace{43.0pt} b \cy_{2,t-1} + \veps_{2t} \label{eq:minvar2}
\end{align}
\end{subequations}
with residuals $\veps_{1t},\veps_{2t} \text{ iid }\sim \normal(0,1)$ so that $\Sigma = I$, and $|a|,|b| < 1$ so that the model is stable.
The transfer function is
\begin{equation}
 	H(z) = \begin{bmatrix} (1-az)^{-1} & cz(1-az)^{-1}(1-bz)^{-1} \\ 0 & (1-bz)^{-1} \end{bmatrix} \label{eq:minvartr}
\end{equation}
From \eqref{eq:cpsdi} the CPSD for $\cy_{1t}$ is $|1-az|^{-2}|1-bz|^{-2} \bracs{1+b^2+c^2 -2b\,\re(z)}$ which factorizes as $H^\rrop(z) \Sigma^\rrop {H^{\rrop *}(z)}$ with $H^\rrop(z)$ of the form $(1-az)^{-1} (1-bz)^{-1} (1-hz)$, so that $\cy_{1t}$ is ARMA($2,1$). We may calculate
\begin{equation}
 	\Sigma^\rrop = \tfrac12 \big(\Delta+ \sqrt{\Delta^2 - 4 b^2}\big) \,, \quad\text{where } \Delta \equiv 1+b^2+c^2 \label{eq:minvargc}
\end{equation}
We have $\gc{\cy_2}{\cy_1} = \ln \Sigma^\rrop$, while causality in the $\cy_1 \to \cy_2$ direction clearly vanishes.

We filter the process \eqref{eq:mvar} through a binomially weighted moving average filter with transfer function
\begin{equation}
 	G(z) = \begin{bmatrix} (1+f_1 z)^r & 0 \\ 0 &  (1+f_2 z)^r \end{bmatrix} \label{eq:wmafilt}
\end{equation}
The filter \eqref{eq:wmafilt} is always causal and stable, and we take $\max(|f_1|, |f_2|) < 1$ so that it is also minimum phase. The filtered series $\by_t \equiv G(z) \cdot \cby_t$ is then ARMA($r,1$) and causalities $y_2 \to y_1$ and $y_1 \to y_2$ are invariant \citep{Barnett:gcfilt:2011}\footnote{While \citet{Barnett:gcfilt:2011} state that causality and stability are sufficient for GC invariance, minimum phase is also required. We thank Victor Solo (private communication) for pointing this out.}.

For our simulation experiment, we set reference parameters
\begin{equation}
	a = 0.9, \quad b = 0.8, \quad c = \sqrt{e^{-F}(e^F-1)(e^F-b^2)} \label{eq:minvar_refparms}
\end{equation}
(so that $\gc{\cy_2}{\cy_1} = F$) with $F = 0.02$ for a causal and $F = 0$ for a null (non-causal) model. Reference parameters for the filters \eqref{eq:wmafilt} were
\begin{equation}
	f_1 = 0.6, \quad f_2 = 0.7 \label{eq:wmafilt_refparms}
\end{equation}
The MA order $r$ is allowed to vary.

\begin{figure}
\includegraphics{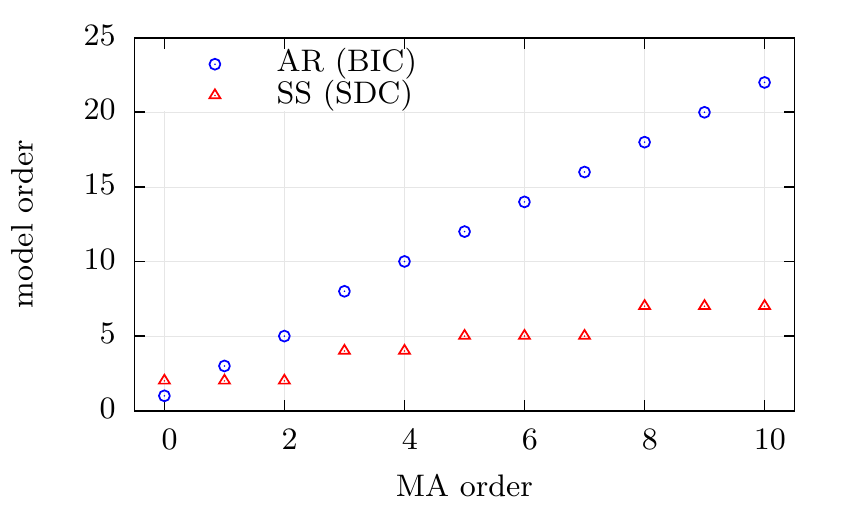}
\caption{AR and SS median model orders plotted against MA order, for the process $\by_t$ with reference parameters \eqreff{eq:minvar_refparms}{eq:wmafilt_refparms}.} \label{fig:minvar_mo}
\end{figure}

To test for power and bias, we obtained empirical distributions for AR and SS estimators $\egc{y_2}{y_1}$, for both null and causal models. Empirical distributions were based on $10,000$ sample simulations of \eqref{eq:mvar} of $T = 1000$ time steps, filtered according to \eqref{eq:wmafilt}, for $r = 0$ (no MA component) to $r = 10$ (strong MA component). For the AR GC estimates, AR modeling of the filtered time series was by OLS. For each sample a model order was estimated using the BIC and $\egc{y_2}{y_1}$ computed from the corresponding autocovariance sequence \citep{Barnett:mvgc:2014}. For SS modeling, the Canonical Correlations Analysis (CCA) SS-SS algorithm \citep{Larimore:1983,Bauer:2005} was used; see \apxref{sec:simulation} for details of the model order selection procedure. SS GC was then calculated according to \eqref{eq:cgctd} using a $\Sigma^\rrop$ estimate obtained by solution of the corresponding DARE as described in the previous section. In both the AR and SS cases, model orders were virtually identical for the null and causal models. \figref{fig:minvar_mo} confirms that AR model order increases (apparently linearly) more rapidly with increasing MA order than SS model order.

\begin{figure*}
\includegraphics{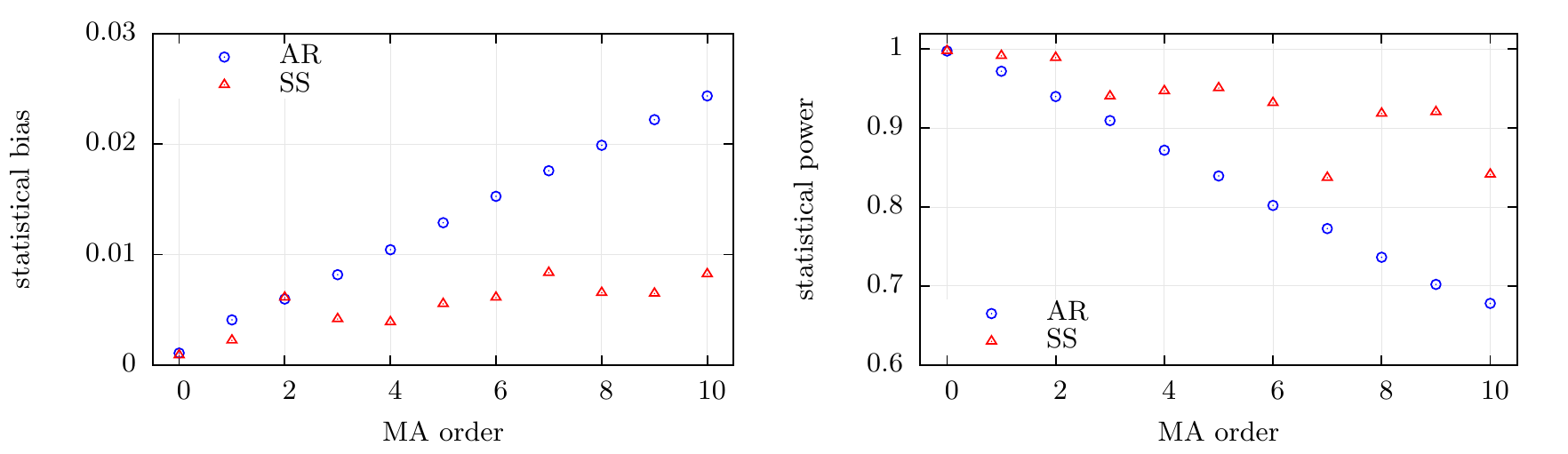}
\caption{Statistical bias (left column), and power at $\alpha = 0.05$ (right column) of $\egc{y_2}{y_1}$ for the process $\by_t$ with reference parameters \eqreff{eq:minvar_refparms}{eq:wmafilt_refparms}, plotted against MA order, as calculated from $\Gamma$ distributions fitted to the empirical distributions.} \label{fig:minvar_gcstats}
\end{figure*}

Empirical null and causal GC distributions were fitted (based on ML parameter estimation) to $\Gamma$ distributions (see \apxref{sec:statinf}). Bias, and statistical power at $\alpha = 0.05$, were calculated from the empirical distributions and also using the fitted $\Gamma$ CDFs; results for the $\Gamma$ fit (\figref{fig:minvar_gcstats}) were very similar to the empirical results. The key outcome is that both bias and statistical power scale more slowly with increasing strength of the MA order for the SS than the AR causal estimators.
We remark that the fluctuations in SS bias and power in \figref{fig:minvar_gcstats} are not statistical fluctuations - nor is there reason to believe that they are artifacts of the experimental procedure (it is, however, not clear whether the SS order estimates are on average truly minimal - see \eg\ the MA order $7$ and $10$ points in \figref{fig:minvar_gcstats}).

\section{Discussion}

In this paper we have described and experimentally validated a new approach to Granger-causal inference, based on state space modeling. This approach offers improved statistical power and reduced bias, when compared to standard autoregressive methods, especially when the data manifest a moving average component. Our primary theoretical contribution is to demonstrate that GC for a state space process, conditional and unconditional, in time and frequency domains, may be simply expressed in terms of SS model parameters, the calculation involving solution of a single DARE. The resulting estimation of GC from empirical time series data is simpler and more computationally efficient and numerically stable than AR-derived causal estimation. Detailed simulation experiments verified the gain in statistical power and reduction in bias resulting from SS as compared to AR modeling of an ARMA process.

These results are important because many application domains furnish time series data likely to feature a moving average component, either intrinsic to the data generating process, or arising from common data manipulation operations such as subsampling, aggregation, additive observational noise, filtering and subprocess extraction. One important example, where the use of GC has remained controversial, is in functional magnetic resonance imaging (fMRI) of the brain, where the time series are highly downsampled and filtered (by hemodynamic processes) reflections of the neural dynamics of interest \citep{ValdesSosaEtal:2011,Handwerker:2012,SethJNEURO:2015}. Other important scenarios where MA components are to be expected include climate science and economics. SS-based Granger causal inference is likely to outperform standard AR-based inference in these and other similar situations, since the class of ARMA (equivalently linear SS) models is closed under these operations. Given the ready availability of easily implementable and computationally efficient SS subspace algorithms, the state space approach should be the default methodology when implementing GC analysis.

The state space approach to Granger-causal inference also seems amenable to extensions beyond linear, stationary and homoscedastic models. For example, subspace algorithms may be adapted for \emph{cointegrated} processes \citep{BauerWagner:2002} and SS models have been developed for heteroscedastic noise variance \citep{WongEtal:2006}. A powerful feature of the Kalman filter underlying state prediction is that it applies to \emph{non}stationary systems, which raises the possibility of a systematic approach to GC analysis for systems where parameters vary over time. Finally, nonlinear state space systems already constitutes a mature field with a large literature. Bilinear SS models in particular form the basis for Direct Causal Modeling (DCM) \citep{FristonEtal:2003}, an approach to causal analysis of coupled dynamical systems sometimes seen as a rival to Granger-causal analysis. The state space setting provides a common framework in which to compare GC and DCM, helping to systematize statistical approaches to characterization of causal connectivity \citep{FristonEtal:2013}.

\subsection*{Acknowledgments}

The authors are grateful to the Dr. Mortimer and Theresa Sackler Foundation, which supports the Sackler Centre for Consciousness Science.

\appendix

\section{GC estimation and statistical inference for AR and SS models} \label{sec:statinf}

A preliminary stage in model identification is choosing an appropriate model order which avoids under- or over-fitting the data. The number of free parameters for an $n$-variable AR model of order $p$ is $p n^2$, while for an SS model with state dimension $m$ it is $2m n$ \citep{HandD:2012}. Standard likelihood-based methods such as the Akaike (AIC) or Bayesian (BIC) information critera \citep{McQuarrie:1998} are commonly employed for model order selection. For both AR and SS models, standard expressions are known for the likelihood function \citep{Lutkepohl:2005}; in the SS case the likelihood is calculated by the \emph{Kalman filter} \citep{AandM:1979}. SS-SS algorithms naturally compute Kalman filter-like states, and also facilitate ``inline'' model selection criteria based on a singular value decomposition (SVD) of a certain (suitably weighted) regression matrix $\bbeta$ \citep{Bauer:2001,GarciaHiernaux:2012}. These criteria may be derived in a single step as part of the parameter estimation procedure, as opposed to ML methods which will generally require multiple parameter estimations over a range of model orders.

In the AR case, a na\"ive estimator for $\cgc{\by_2}{\by_1}{\by_3}$ can be formed by performing separate regressions for the full and reduced AR models to obtain estimates of the full and reduced covariance matrices $\Sigma$ and $\Sigma^\rrop$ respectively, which appear in the expression \eqref{eq:cgctd} for (time domain) GC. However, this is not a recommended procedure \citep{Dhamala:2008a,Barnett:mvgc:2014}; for conditional spectral Granger causality estimation in particular, independent full and reduced parameter estimates may lead to spurious, or even negative causality estimates \citep{Chen:2006} \footnote{We do not recommend the ``partition matrix'' technique proposed in \cite{Chen:2006} to address this issue, as the algorithm is technically flawed and leads to inaccurate results. The flaw relates to an inappropriate use of filtering in place of spectral factorization.}. This problem relates to the fact that the class of finite order AR models is not closed under subprocess extraction; even if the full AR model is of finite order, the reduced model will generally be ARMA.

Alternative approaches are to derive reduced model parameters from the CPSD \citep{Dhamala:2008a,Dhamala:2008b} or autocovariance sequence \citep{Barnett:mvgc:2014}, either of which may be calculated directly from the full AR parameter estimates. Although these single-regression methods generally have greater statistical power, the computations involve various approximations impacting numerical accuracy and stability, and are potentially computationally expensive. Computation of Granger causalities from the autocovariance sequence involves approximation of the autocovariance to a potentially large number of lags (specifically if autocovariance of the estimated AR model decays slowly) and application of Whittle's recursive time-domain spectral factorization algorithm \citep{Whittle:1963}. Whittle's algorithm, while generally stable and accurate, does not scale particularly well with system size. Computation of causalities from the CPSD involves the choice of a potentially high frequency resolution (again depending on autocovariance decay) and spectral factorization via Wilson's iterative frequency-domain algorithm \citep{Wilson:1972}. The latter again scales poorly with system size and, although in theory the iteration converges quadratically, may, in the authors' experience, suffer from stability and numerical accuracy issues. See \citet{Barnett:mvgc:2014} for further discussion on these issues.

In the SS case, Granger causalities could in principle be computed similarly from the autocovariance sequence or CPSD. We may calculate that for the innovations form SS model \eqref{eq:ssi}, the autocovariance sequence $\Gamma_k \equiv \expect{\by^{\phantom\trop}_t \by^\trop_{t-k}}$ is given by
\begin{subequations}
\begin{align}
	\Gamma_0 &= C \Omega C^\trop + \Sigma \label{eq:acpo0} \\
	\Gamma_k &= C A^k \Omega C^\trop + C A^{k-1} K \Sigma, \qquad k > 0 \label{eq:acpok}
\end{align} \label{eq:acpo}%
\end{subequations}
where $\Omega \equiv \expect{\bx^{\phantom\trop}_t \bx^\trop_t}$ satisfies the discrete Lyapunov equation $\Omega = A \Omega \tA + K \Sigma K^\trop$, while the CPSD is given by \eqreff{eq:ssioar}{eq:cpsdi}. However, as remarked above, there are drawbacks to these approaches, which are in any case unnecessary, since we have shown how reduced SS model parameters, and thence Granger causalities, may be calculated directly from full SS model parameters. This computation, as we have seen, involves the solution of a single (or two, if the original SS model is not already in innovations form) DARE, for which efficient, stable and accurate algorithms exist (\eg\ \citep{ArnoldLaub:1984}), and is thus likely to be less computationally intensive as well as more stable and numerically accurate than CPSD- or autocovariance-based methods. The efficacy of this method is demonstrated in the simulation experiments.

\begin{figure*}
\includegraphics{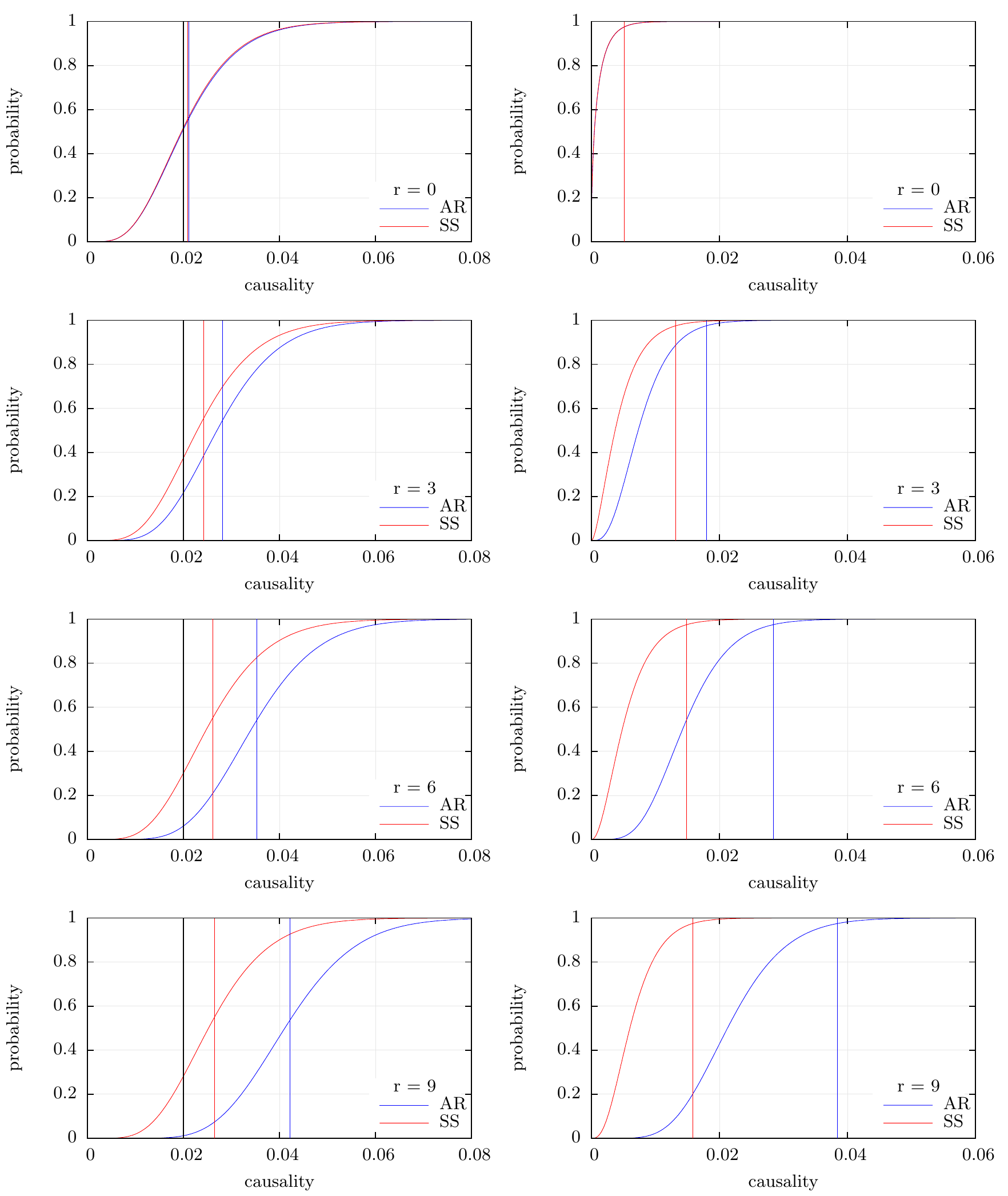}
\caption{Empirical distribution (CDF) of $\egc{y_2}{y_1}$ for the AR($1$) process \eqref{eq:mvar} with MA filter \eqref{eq:wmafilt}, with reference parameters \eqreff{eq:minvar_refparms}{eq:wmafilt_refparms} and $T = 1000$ time steps. Left column: causal model ($F = 0.02$, black vertical line is actual causality, colored vertical lines indicate means), right column: null model ($F = 0$, colored vertical lines indicate critical causalities $F_{crit}$ at $95\%$ significance). MA order $r$ increases from top to bottom.} \label{fig:minvar_gcdist}
\end{figure*}

Regarding statistical inference, if in the AR case the full and reduced models are independently estimated then the associated estimator of $\cgc{\by_2}{\by_1}{\by_3}$ is the log-likelihood ratio statistic for a nested AR model under the null hypothesis of vanishing causality. As such, from the standard large sample theory \citep{Wilks:1938,Wald:1943}, it has an asymptotic sampling distribution under the null hypothesis as a $\chi^2$ with $d_f = p n_1 n_2$ degrees of freedom, where $n_1,n_2$ are the dimensions of $\by_1,\by_2$ respectively and $p$ the AR model order. For $\cgc{\by_2}{\by_1}{\by_3}$ estimated from the autocovariance sequence, we have found that the estimator is well-approximated by a $\Gamma$ distribution (of which the $\chi^2$ distribution is a special case) with shape parameter corresponding to (possibly non-integer) $d_f < pn_1n_2$. See \figref{fig:minvar_gcdist},
where empirical null and causal AR GC distributions for our simulation experiment were fitted, based on ML parameter estimation, to $\Gamma$ distributions; the fitted CDFs are indistinguishable from the empirical CDFs at the scale of the figure. We have not been able to ascertain a simple formula for $d_f$ in terms of model parameters.

In the SS case, since the full and reduced SS models are not simply nested, it is not clear how the large sample theory might obtain \footnote{Although the non-causality condition \eqref{eq:cgcnull} suggests $d_f = m n_1 n_2$ degrees of freedom for an appropriate nested model, preliminary testing with SS-SS GC estimates indicates that, even for quite long time series, the corresponding $\chi^2$ distributions do not in general furnish a good fit with the empirical distributions.}. Under the null hypothesis of zero causality, the sampling distribution for an estimator of $\cgc{\by_2}{\by_1}{\by_3}$ obtained by replacing $\Sigma_{11}$ and $\Sigma^\rrop_{11}$ (the latter derived as described in \secref{sec:gc}) model estimates, is again well-approximated by a $\Gamma$ distribution. In \figref{fig:minvar_gcdist} fitted CDFs are again indistinguishable from the empirical CDFs.

Lacking known theoretic asymptotic distributions, permutation and bootstrapping (or other surrogate data techniques) may be deployed for significance testing and estimation of confidence intervals respectively. As regards frequency-domain Granger causality estimators, as far as we are aware no sampling distributions are known, even in the AR case.

\section{Statistical power and bias for GC estimators} \label{sec:statpow}

Statistical power of an estimator is defined as the fraction of \emph{true positives} (\ie\ correct rejections of the null hypothesis) obtained at a given significance level. Suppose that the cumulative distribution function (CDF) of a GC estimator $\EGC$ given an \emph{actual} causality $F$ is $\Phi_F(\ldots)$. To test for statistical significance at level $\alpha$, we would reject the null hypothesis of zero causality when $\EGC > F_{crit} \equiv \Phi_0^{-1}(1-\alpha)$. Given an actual causality $F$, then, the statistical power is $\prob{\EGC > F_{crit}} = 1-\Phi_F\big(\Phi_0^{-1}(1-\alpha)\big)$. Under the same assumptions, the bias of the estimator is defined as $\expect{\EGC}-F$.

\section{SS-SS modeling and model order selection} \label{sec:simulation}

For SS GC estimates, SS modeling was performed using the CCA state space-subspace algorithm \citep{Larimore:1983,Bauer:2005}. We note firstly that not much appears to be known about optimal values for the past and future horizons for subspace algorithms. \citet{BauerWagner:2002} suggest setting both to $2p_{\text{AIC}}$ where $p_{\text{AIC}}$ is the AIC-based AR model order estimate (see also \citep{GarciaHiernaux:2012}). In this study we used $p_{\text{BIC}}$, the BIC AR order estimate, which was found to yield better results.

\begin{figure}
\includegraphics{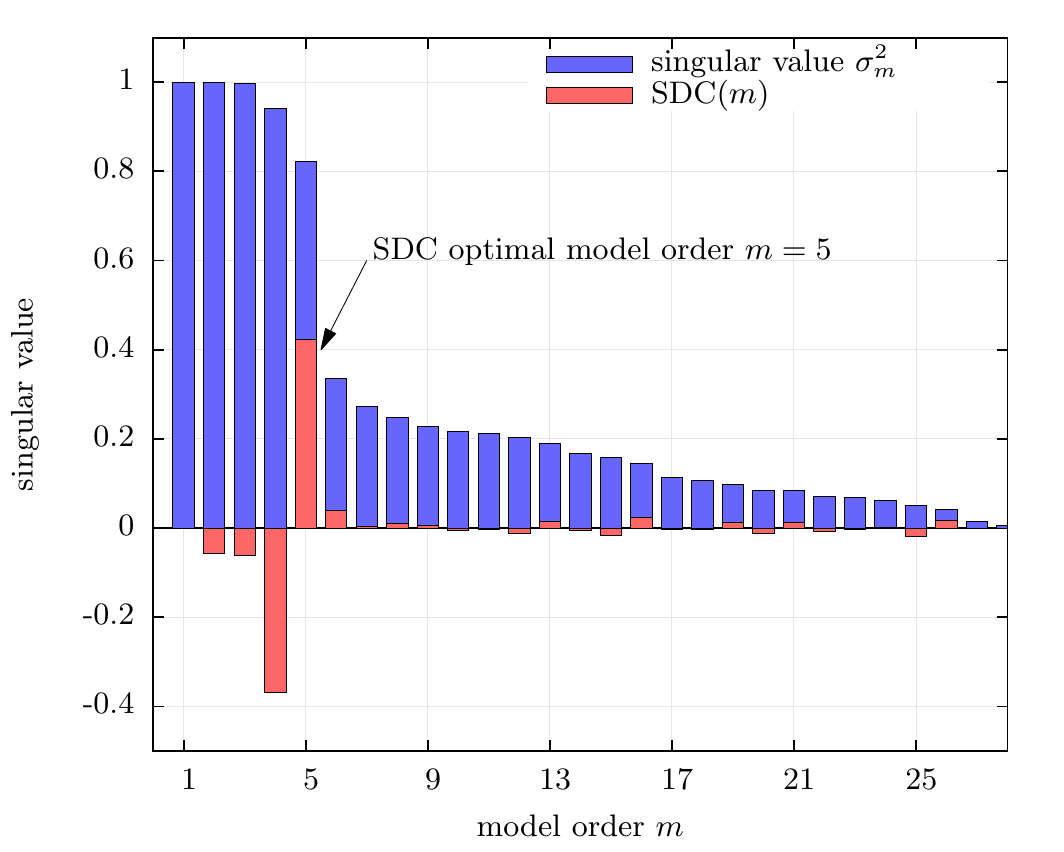}
\caption{The SS model order criterion SDC$(m)$ for a sample filtered process $\by_t$ (see text for details).} \label{fig:minvar_SVC}
\end{figure}

To estimate model order, we used a ``steepest drop criterion'' (SDC), defined as $\text{SDC}(m) \equiv \sigma_{m+1}^2 -2 \sigma_m^2 + \sigma_{m-1}^2$, where $\sigma_m^2$ is the $m$th singular value of the CCA-weighted regression matrix $\bbeta$ calculated in the SS-SS algorithm \citep{Bauer:2001}. The optimal model order estimate is then taken as the value of $m$ which maximizes SDC$(m)$. The criterion thus attempts to pinpoint the steepest drop in singular value with increasing state space dimension (\figref{fig:minvar_SVC}). We also tested alternative model order selection schemes, including the NIC, SVC and IVC of \citet{Bauer:2001}. SDC estimates were, overall, similar to those of SVC, but proved more robust to variation between sample time series.

\bibliography{ssgc}

\end{document}